\documentclass[runningheads]{llncs}
\usepackage{amsfonts, amsmath, amstext, amssymb}

\usepackage[dvips]{graphicx}
\usepackage{latexsym}
\usepackage{pgf,pgfarrows,pgfnodes,pgfautomata,pgfheaps}
\usepackage{url}

\usepackage{color}

\newtheorem{thm}{Theorem}[section]
\newtheorem{defn}[thm]{Definition}
\newtheorem{prop}[thm]{Proposition}

\newtheorem{cor}[thm]{Corollary}

\def\bec{\begin{cor}}
\def\enc{\end{cor}}

\def\bet{\begin{thm}}
\def\ent{\end{thm}}

\def\becor{\begin{cor}}
\def\encor{\end{cor}}

\def\bel{\begin{lem}}
\def\enl{\end{lem}}

\def\bedef{\begin{defn}}
\def\endef{\end{defn}}

\def\bep{\begin{prop}}
\def\enp{\end{prop}}

\newenvironment{pf}{\begin{trivlist}\item[\hskip\labelsep
{\it Proof.}]}{\end{trivlist}}

\newcommand{\set}[2]{\ensuremath{ \{ #1 : #2 \} }}

\newcommand{\DS}[2]{\ensuremath{\text{DgSp}_{#1}(#2)}}
\renewcommand{\deg}[1]{\ensuremath{\text{deg}(#1)}}

\newcommand{\Z}{\mathbb{Z}}

\newcommand{\Q}{\mathbb{Q}}

\newcommand{\C}{\mathcal{C}}
\newcommand{\A}{\mathcal{A}}
\newcommand{\B}{\mathcal{B}}

\newcommand{\Cbar}{\overline{C}}
\newcommand{\Dbar}{\overline{D}}

\newcommand{\Qbar}{\overline{\mathbb{Q}}}

\newcommand{\Gal}[2]{\text{Gal}(#1/#2)}
\newcommand{\la}{\langle}
\newcommand{\ra}{\rangle}

\newcommand{\comment}[1]{}

\def\converges{\!\downarrow}

\newcommand{\dom}[1]{\text{dom}(#1)}
\newcommand{\rg}[1]{\text{range}(#1)}

\newcommand{\bfc}{\boldsymbol{c}}
\newcommand{\bfd}{\boldsymbol{d}}
\def\bfz{\boldsymbol{0}}
\def\s01{\ensuremath{\Sigma^0_1}}
\def\d02{\ensuremath{\Delta^0_2}}
\def\phi{\varphi}

\def\res{\!\!\upharpoonright\!}

\begin{document}

\title{Degree Spectra for Transcendence in Fields\thanks{The research
of the first author is supported by RSF Grant no. 18-11-00028; he is also
funded by the Russian Ministry of Education and
Science (project 1.451.2016/1.4) as a federal professor in mathematics.
The second author was partially supported by Grant \# 581896 from
the Simons Foundation, and the second and third authors were both supported
by grants from the City University of New York PSC-CUNY Research Award Program.
The authors wish to acknowledge useful conversations with Dr.\ Kenneth Kramer.
}}
\author{Iskander Kalimullin\inst{1}
\and Russell Miller\inst{2,4}
\and Hans Schoutens\inst{3,4}
}
\authorrunning{Iskander Kalimullin, Russell Miller, \& Hans Schoutens}
\institute{Kazan Federal University, Kremlyovskaya St.\ 18, 420008 Kazan, Russia\\
\email{ikalimul@gmail.com}
\and Queens College, 65-30 Kissena Blvd., Queens, NY 11367, U.S.A.\\
\email{Russell.Miller@qc.cuny.edu}
\and New York City College of Technology, 300 Jay Street, Brooklyn, NY 11201 U.S.A.\\
\email{hschoutens@citytech.cuny.edu}
\and C.U.N.Y.\ Graduate Center, 365 Fifth Avenue, New York, NY  10016 U.S.A.
}

\maketitle

\begin{abstract}
We show that for both the unary relation of transcendence
and the finitary relation of algebraic independence on a field,
the degree spectra of these relations may consist of any
single computably enumerable Turing degree, or of
those c.e.\ degrees above an arbitrary
fixed $\Delta^0_2$ degree.  In other cases, these spectra
may be characterized by the ability to enumerate an arbitrary $\Sigma^0_2$ set.
This is the first proof that a computable field can fail to have
a computable copy with a computable transcendence basis.

\end{abstract}

\section{Introduction}
\label{sec:intro}

It has been known since the work of Metakides and Nerode in
\cite{MN79} that a computable field need not have a computable
transcendence basis.  This result, readily established, is fundamental
to the study of effectiveness for fields.  Under the usual definition
of computable structure, a computable field is simply a field
whose domain is a computable subset of $\omega$ (usually just
$\omega$ itself, the set of all nonnegative integers) and whose atomic
diagram, in the language with addition and multiplication, is computable.
The theorem of Metakides and Nerode shows that working with
an arbitrary computable field will be difficult, as one cannot in
general distinguish the algebraic elements of the field (relative
to the prime subfield, either $\Q$ or $\Z/(p)$) from those
transcendental over the prime subfield.

We recently realized that the following very natural question
had not been addressed:  must every computable field be
isomorphic to a computable field with a computable transcendence basis?
It is well known that there need not exist a computable isomorphism
between two isomorphic computable fields,
and so it is plausible that the answer might be
affirmative:  two computable fields, one with a computable
transcendence basis and the other without any such basis,
can certainly be isomorphic.  (We normally refer to isomorphic
computable fields as \emph{computable copies} of each other.)
In case of an affirmative answer, one would be justified in always
assuming a computable transcendence basis, as this would only
require choosing a ``nice'' computable copy of the field in question.

Initially we were optimistic that the answer would indeed be affirmative,
and even that a single Turing procedure might produce such a copy
uniformly, using the original field's atomic diagram as an oracle.  (This
could also then be extended to noncomputable fields.)  However,
conversations with our colleague Ken Kramer disabused us of that
notion, and in fact we will demonstrate here that the answer is
negative.  For uniform procedures, the negative answer is proven
in Section \ref{sec:genus}, which introduces and illustrates the use
of algebraic curves of positive genus for this purpose.  The remainder
of the article shows that there is not even any nonuniform procedure:
certain computable fields have no computable copy with any
computable transcendence basis.  Indeed, the spectrum of the transcendence
relation on a field has many possible configurations, plenty of which
do not include the degree $\bfz$.  For each computably enumerable
Turing degree $\bfc$, it is possible for transcendence to be
intrinsically of degree $\bfc$, or for it intrinsically to compute $\bfc$,
or for it intrinsically to enumerate a given $\Sigma^0_2$ set.
The proofs here make substantial use of results on algebraic curves
developed during earlier work by two of us in \cite{MS13}.

Transcendence bases are not in general definable, and a single field
of infinite transcendence degree will have continuum-many different
transcendence bases.  This makes it difficult to define ``the'' Turing degree
for transcendence bases.  To address this, we use two $L_{\omega_1\omega}$-definable
relations on fields.  The transcendence relation $T$, which is unary, holds of those elements
not algebraic over the prime subfield $Q$ of a field $F$:
$$ x\in T \iff (\forall f\in Q[X]^*) f(x)\neq 0.$$
(Here $Q[X]^*$ is the set of nonzero polynomials over $Q$.)
The algebraic independence relation $I$ is a generalization of this to tuples
of all arities $n$:
$$ (x_1,\ldots,x_n)\in I \iff (\forall f\in Q[X_1,\ldots,X_n]^*) f(x_1,\ldots,x_n)\neq 0.$$
A computable field of infinite transcendence degree will possess transcendence bases
having each Turing degree $\geq_T\deg{I}$, but not of any other Turing degree:
given any basis as an oracle, one can decide the independence relation on the field,
and conversely, from an $I$-oracle, one can compute a transcendence basis for $F$.
Thus the Turing degree of $I$ can stand in for the set of Turing degrees of transcendence
bases, as this set is the upper cone above $\deg{I}$.  In turn, the unary relation $T$
is always computable from $I$, although sometimes strictly below $I$ under Turing reducibility.
In the fields we consider here, we will always have $I\equiv_T T$.  We remark the
following useful property.
\begin{lemma}
\label{lemma:cebasis}
In a computable field $F$, for a Turing degree $\bfd$, every $\bfd$-computably
enumerable transcendence basis $B$ over the prime subfield $Q$ is $\bfd$-computable.
\end{lemma}
\begin{pf}
Given any $x\in F$, use a $\bfd$-oracle to list the elements $b_0,b_1,\ldots$ of $B$,
and search for an $n$ and a polynomial $f\in Q[X,Y_0,\ldots,Y_n]^*$
with $f(x,b_0,b_1,\ldots,b_n)=0$ in $F$.  This search must terminate,
and $x\in B$ just if $x\in\{b_0,\ldots,b_n\}$.
\qed\end{pf}

Our notation is generally standard.  The articles \cite{MNotices,M11} form good introductions
to computable field theory, and myriad other articles have contributed to the area:
\cite{E77,FS56,MN79,R60,vdW70} all have historical importance, while
\cite{FKM09,MPSS18,MS13} describe related recent work in the discipline.

\section{Curves of Positive Genus}
\label{sec:genus}

\begin{proposition}
\label{prop:genus}
Let $f(Y,Z)=0$ define a curve over a field $k$ of characteristic $0$.
If the genus of this curve is positive, then $f=0$
has no solutions in the purely transcendental field extension
$K=k(t_1,t_2,\ldots)$ except those solutions within $k$ itself.
(We say that $f=0$ has no \emph{nonconstant} solutions in $K$.)
\end{proposition}
\begin{pf}
A solution in $K$ would lie within some subfield $k(t_1,\ldots,t_n)$,
so we prove by induction on $n$ that no such subfield contains a nonconstant solution.
For $n=1$, we can view the extension $k(t_1)$ as an algebraic curve
of genus $0$ over $k$.  If an extension $k(y,z)$ (where $f(y,z)=0$)
lies within $k(t_1)$, with $\{ y,z\}\not\subseteq k$, then the
Riemann-Hurwitz formula dictates that $f=0$ must also have genus $0$,
contradicting the hypothesis of the proposition.

For the inductive step, suppose $y,z\in k(t_1,\ldots,t_{n+1})$ satisfy
$f(y,z)=0$.  We express $y=\frac{g_1}{h_1}$ and $z=\frac{g_2}{h_2}$
as rational functions of $t_1,\ldots,t_{n}$ over the field $k(t_{n+1})$
Of course, each of $g_1,g_2,h_1,h_2$ has finitely many nonzero coefficients
in that field, and the pairs $(g_1,h_1)$ and $(g_2,h_2)$ may be taken to have
no common factor.  Having characteristic $0$, $k$ is infinite, so
it must contain an element $a$ such that, when $t_{n+1}$ is replaced by $a$,
all of these coefficients remain nonzero and no common factors are introduced.
Substituting $a$ for $t_{n+1}$ in $y$ and $z$ yields a solution to $f=0$
in $k(t_1,\ldots,t_n)$.  By inductive hypothesis this solution lies in $k$,
meaning that the original $y$ and $z$ did not involve any of $t_1,\ldots,t_n$
(since no common factors arose to be cancelled when we made the substitution).
But then $y$ and $z$ were a solution to $f=0$ in $k(t_{n+1})$, hence
must be a constant solution, according to the base case of the induction.
\qed\end{pf}

\begin{corollary}
\label{cor:genus}
There is no uniform procedure for transforming a countable field into
an isomorphic countable field that decides its own transcendence relation $T$.

That is, there does not exist any Turing functional $\Psi$ such that, for every
atomic diagram $F$ of a countable field with domain $\omega$ (in the signature
with just $+$ and $\cdot$), $\Psi^F$ computes the atomic diagram,
in the larger signature with $+$, $\cdot$ and $T$, of a structure with reduct $F$
and for which $T$ is the (unary) transcendence relation.  The same holds with
the ($<\!\omega$-ary) algebraic independence relation $I$ in place of $T$.
\end{corollary}
\begin{pf}
Suppose $\Psi$ were such a functional.  Fix an irreducible curve with affine equation $f(Y,Z)=0$
of positive genus over $\Q$, and let $F$ be a presentation of the field
$\Q(y_1,y_2,y_3,\ldots)(z_1)$, with $\{y_1,y_2,\ldots\}$ algebraically
independent over $\Q$ and $f(y_1,z_1)=0$.  Then $\Psi^F$ must compute
the atomic diagram of an isomorphic field $L$ with a transcendence relation $T$,
say with isomorphism $h:F \to L$.
Therefore $T$ will hold of the $1$-tuple $h(y_1)$ in the field $L$.
Let $\sigma$ be an initial segment of the atomic diagram of $F$ such that
$\Psi^{\sigma}$ ensures that $T$ holds of $h(y_1)$ and that $f(h(y_1),h(z_1))=0$ in $L$.

Now let $E$ be a presentation of the field whose atomic diagram begins with $\sigma$.
However, the atomic diagram of $E$ (beyond $\sigma$) specifies that $y_1$
is in fact rational itself, in some way consistent with $\sigma$, and thus that $z_1$
is algebraic over $\Q$ in $E$.  (It may not be possible to make $z_1$ rational too,
as $f=0$ will have only finitely many solutions in $\Q$, by Faltings' proof of the
Mordell Conjecture.  However, there is no difficulty in making $z_1$ algebraic over $\Q$.)
The rest of $E$ is then generated by this portion and by elements $y_2,y_3,\ldots$
algebraically independent over $\Q$, just as in $F$.

Thus $\Psi^E$ will build a field in which the domain element $h(y_1)$ is transcendental
over $\Q$ (being thus specified by $\Psi^E$) and $f(h(y_1),h(z_1))=0$.
However, $E$ is a purely transcendental extension of the field $k=\Q(z_1)$,
which is algebraic over $\Q$.  By Proposition \ref{prop:genus}, $E$ does not contain
any solution to $f(Y,Z)=0$ outside of $k$, so every solution in $E$ consists of elements
algebraic over $\Q$.  This ensures that $E$ and the field with atomic
diagram $\Psi^E$ are not isomorphic as fields, proving the Corollary.
(The result for the relation $I$ follows directly.)
\qed\end{pf}

\section{Background on Algebraic Curves}
\label{sec:Fermat}

Corollary \ref{cor:genus} proved that there is no uniform method of
taking a computable field and producing a
computable copy with a computable transcendence basis.  Now we wish
to show that a single computable field can entirely fail to have a computable
copy with a computable transcendence basis.  Indeed, we will establish far
more specific results, with detailed descriptions of the possible degrees
of transcendence bases in computable copies of the field.  To do this, however,
we need to work with infinitely many curves of positive genus at once, as a single
curve will only allow our field to avoid being isomorphic to a single
computable field with computable transcendence basis.

Fortunately, an appropriate collection of curves has already been built.
We recall the following result from \cite{MS13}, as stated there.
\begin{theorem}[Miller \& Schoutens, Thm.\ 3.1 of \cite{MS13}]
\label{thm:noncovering}
There exists a non-covering collection $\C$ of curves with the effective Faltings property,
over a computable ground field $k$.
\end{theorem}
That is, $\C=\{ f_0,f_1,\ldots\}$ is an infinite set of smooth projective curves $C_i$
with corresponding affine equations $f_i\in k[X,Y]$ such that:
\begin{itemize}
\item
for each $i$, the function field $k(f_i)$ does not embed into $k(\C-\{ f_i\})$; and
\item
the function $i\mapsto f_i(k)$ giving (a strong index for the finite set of) all
solutions of $f_i(X,Y)=0$ in $k^2$ is computable.  (This is the \emph{effective
Faltings property}.)
\end{itemize}
In particular, the second item requires that each equation $f_i=0$
should have only finitely many solutions by elements of $k$.

The specific example $\C$ given in \cite{MS13} is in fact a collection
of Fermat curves $f_i=X^{q_i}+Y^{q_i}-1$, for a fixed increasing computable
sequence $q_0<q_1<\cdots$ of odd prime numbers.  By Fermat's Theorem,
each has exactly two solutions in $k=\Q$, and the non-covering property
for this $\C$ is established in \cite{MS13}.  It is believed that many other
computable sets of curves have the same property, but rather than pursuing
that question here, we will use this same set $\C$.  It should be borne in mind
that not all odd primes belong to the sequence $\la q_i\ra_{i\in\omega}$.
Indeed, this sequence is quite sparse within the primes:  each element
is the least prime $q_{i+1}>(4(q_i-1)(q_i-2))^2$, with $q_0$ chosen to
be $5$ (or any other odd prime except $3$, which is ruled out because
the genus of $(X^d+Y^d-1)$ is $\frac{(d-2)(d-1)}2$ and we need genera $>1$).

\comment{
The existence is such a prime is ensured at each step by Dirichlet's Theorem
on primes in arithmetic sequences.  This choice will allow us to apply the following lemma.

\begin{lemma}
\label{lemma:Galois}
For the sequence $\la q_i\ra$ given above, let $F\supseteq\Q$ be the infinite algebraic extension
generated by the real roots $\sqrt[q_j]{a_j}$ for all $j\neq i$, where each $a_j$
may be any rational number one wishes.  Then $F$ contains no nontrivial solution
to $f_i(X,Y)=0$.
\end{lemma}
\begin{pf}
We view $F$ as the result of an infinite sequences of algebraic extensions,
each one of prime degree $q_j$ (for $j\neq i$) over the previous one.  If $F$
contains a solution $(x,y)$ to $f_i$, then $x$ and $y$ must appear in the same extension,
for if $x$ were in one extension and $y$ in a later one, then $y$ would have
had degree $q_j$ over $x$, and $j\neq i$.  So there is an extension $K=L(y_j)\supseteq L$
within this sequence of some degree $q_j$, with $x$ and $y$ both in $(K-L)$.

Let $N$ be the smallest normal extension of $L$ containing $K$.  Then $N$
is generated over $L$ by all the $q_j$-th roots of some rational number $r=1-x_j^{q_j}$,
and so $[N:L]=q_j(q_j-1)$.  By our choice of $q_j$ (if $i<j$) or of $q_i$ (if $j<i$),
$q_i$ cannot divide $[N:L]$.  However, by normality, $N$ contains all the
$q_i$-th roots of $(1-x^{q_i})$, including $y$, and $\Gal(N/L)$ acts transitively
on these roots, forcing $q_i$ to divide $[N:L]$, which yields a contradiction.
\qed\end{pf}
}

For the rest of this article we fix these curves $C_0,C_1,\dots$ with affine equations
$f_0,f_1,\ldots$ exactly as given here.
The usefulness of Theorem \ref{thm:noncovering} lies in the fact that it enables
us to adjoin to a ground field $k$ (such as $\Q$ or $\Qbar$)
a transcendental element $x$ and then an element
$y$ satisfying $f_i(x,y)=0$ (so $y$ is also transcendental) without creating
any transcendental solutions to any other $f_j$ in the new field $K=k(x)[y]/(f_i)$.
Indeed, our $k$ might already have been built this way, with pairs $(x_0,y_0),
\ldots, (x_{i-1},y_{i-1})$ of transcendental solutions to $f_0,\ldots,f_{i-1}$, say,
and the new $(x,y)$ will not generate any solutions to any of $f_0,\ldots,f_{i-1}$
that were not already in $k$.  This allows us to work independently with the
distinct polynomials $f_i$ and their solutions, and avoids the need for priority
arguments and the like.

It should be noted that $\Q(x_i)[y_i]/(f_i)$ actually contains eight solutions to $f_i(X,Y)=0$.
Two are the trivial solutions $(0,1)$ and $(1,0)$, which we can always recognize and ignore.
Then the solution $(x_i,y_i)$ generates $(\frac{-x_i}{y_i},\frac1{y_i})$, $(\frac{-y_i}{x_i},\frac1{x_i})$,
and the transpositions of these three.  Moreover, we will sometimes work over
the algebraic closure $\Qbar$ instead of $\Q$, and in this case we get $6q_i^2$
nontrivial solutions from $(x_i,y_i)$, since for a primitive $q_i$-th root $\theta$
of unity, $(x_i\theta^j,y_i\theta^k)$ will be another solution.  (There will also
be plenty of non-transcendental solutions within $\Qbar$, of course.)  However,
it was shown by Leopoldt \cite{L96} and Tzermias \cite{T95} that no further transcendental
solutions will exist; this result is also used in \cite{MS13}, appearing there as Theorem 4.4.

\section{Examples of Degree Spectra}
\label{sec:EXs}

Our initial goal was to produce a computable field such that no computable
copy of the field has a computable transcendence basis.  In fact, we will give
a much more specific answer to the question, using the well-established notion
of the degree spectrum of a relation.

\begin{definition}
\label{defn:spec}
For a computable structure $\A$ and a relation $R$ on $\A$, the
\emph{Turing degree spectrum of $R$ on $\A$} is the set
of all Turing degrees of images of $R$ in computable structures isomorphic to $\A$:
$$ \DS{\A}R = \set{\deg{g(R)}}{g:\A\to\B\text{~is an isomorphism onto a computable~}\B}.$$
\end{definition}
In many contexts this definition is restricted to $n$-ary relations $R$,
but it applies equally well to finitary relations, i.e., those defined
on all finite tuples from $(\A)^{<\omega}$, of arbitrary length.

Notice first that both the (unary) transcendence relation $T$ on a computable
field and the (finitary) relation $I$ of algebraic independence are definable
in the field by computable infinitary $\Pi^0_1$ formulas.  Therefore,
in every computable field $F$, both $T$ and $I$ will be $\Pi^0_1$
sets, hence of c.e.\ Turing degree.  This places an upper bound
on the complexity of the degrees in $\DS{F}T$ and $\DS{F}I$,
since all such degrees must be $\leq_T\bfz'$.  Even below that bound,
it also rules out a number of further candidates, namely those degrees
that do not contain any c.e.\ set.  (We call these \emph{properly $\Delta^0_2$
degrees}, meaning that they are $\Delta^0_2$ but not $\Sigma^0_1$.)

Our ultimate goal is to know the degrees of the transcendence bases
in the various computable copies of $F$.
Recall, however, that these are precisely the degrees $\geq_T\deg{I}$.
Thus, once one knows the Turing degree of $I$ in a particular computable
copy of the field, one knows all the degrees of transcendence bases
in that copy, and so we view $\DS{F}I$ as a reasonable answer to the question.

Our first example shows that, for a fixed computably enumerable Turing degree
$\bfc$, the relations $T$ and $I$ can both be \emph{intrinsically of degree $\bfc$}.
This term was used in \cite{DM91}, in which Downey and Moses showed that
the relation of adjacency in a computable linear order can be intrinsically of degree $\bfz'$.
Subsequently, Downey, Lempp, and Wu showed in \cite{DLW10} that
the only degrees $\bfc$ for which the adjacency relation can be intrinsically of degree
$\bfc$ are $\bfc=\bfz'$ and (if the adjacency relation is finite) $\bfc=\bfz$.
Therefore Theorem \ref{thm:singleton} distinguishes the situation for transcendence
and for independence in fields from that of adjacency in linear orders.

\begin{theorem}
\label{thm:singleton}
For each computably enumerable Turing degree $\bfc$, there exists
a computable field $F$ for which the spectrum of the transcendence relation
$T$ and of the independence relation $I$ are both the singleton $\{\bfc\}$.

\end{theorem}
\begin{pf}
Fixing a computable enumeration $\la C_s\ra_{s\in\omega}$ of a c.e.\ set $C\in\bfc$,
we can describe the isomorphism type of our field quickly.  For every $i\notin C$,
it will contain a
transcendental element $x_i$ (over the ground field $\Q$) and an additional element
$y_i$ with $f_i(x_i,y_i)=0$.  Moreover, these elements $x_i$ will form a
transcendence basis, as $i$ ranges over $\Cbar$.  For each $i\in C$,
the field will contain elements called $x_i$ and $y_i$, again satisfying $f_i(x_i,y_i)=0$,
but this $x_i$ will lie within $\Q$, making $y_i$ algebraic over $\Q$.  These
$x_i$ and $y_i$ (for all $i\in\omega$) will generate the field.  (The choice of which
rational number equals $x_i$, for $i\in C$, will depend on the least $s$ with $i\in C_s$.)

Next we give a computable presentation $F$ of this field.
At stage $0$, $F_0$ consists of a finite substructure of the field $\Q$ (with the operations
viewed as relations, so that it makes sense for $\Q$ to have a finite substructure).
At stage $s+1$, we add elements $x_s$ and $y_s$ to $F_s$, along with as many
new elements as are needed in order for the relational atomic diagram of $F_{s+1}$
to specify that $f_s(x_s,y_s)=0$ (but without making $x_s$ itself algebraic over $F_s$).
Then, for the least $i\leq s$ (if any) such that
$i\in C_s$ and we have not yet acted on behalf of $i$, we add enough new
elements to $F_{s+1}$ and define the operations on them to make
$x_i$ lie within $\Q$ (in a way consistent with $F_s$, of course:  nothing in the
atomic diagram should ever be redefined).  This must be possible, since
$x_i$ has been treated as a transcendental up until this stage.
Finally, we take another step to close $F$ under the field operations,
adding another element and extending the relations in $F_{s+1}$
in a way consistent with the principle that the set
$$\set{x_i}{i\leq s~\&~\text{we have not yet acted on behalf of~}i}$$
should form a transcendence basis for $F_{s+1}$. That is, we make sure
not to create any algebraic relations involving these $x_i$'s, and all
new elements added to the field are generated by $\set{x_i,y_i}{i\leq s}$.
This is the entire construction, and it is clear that it does ultimately build
a computable field $F=\cup_s F_s$.  Furthermore, $\set{x_i}{i\notin C}$
will indeed be a transcendence basis for $F$, and every $x_i$ with $i\in C$
will lie within $\Q$ in the field $F$.

We now argue that for every computable field $E\cong F$, the transcendence
relation $T$ and the independence relation $I$ on $E$ have $T\equiv_T I\equiv_T C$.
First, given a $C$-oracle, we enumerate a transcendence basis for $F$ by
collecting, for each $i\notin C$, the first coordinate of the first pair $(x,y)$
that we find in $E^2$ for which $f_i(x,y)=0\neq xy$.  (This pair is not unique,
as mentioned in Section \ref{sec:Fermat}, but the six possible $x$-values are pairwise
interalgebraic.)  By the construction, this enumerates a transcendence basis $B$
for $F$, hence computes one, by Lemma \ref{lemma:cebasis}, and from $B$
we can compute $I$ and $T$.

To show that $C\leq_T T$, we claim that $i\in\Cbar$ just if there
exists a pair $(x,y)$ of elements of $E$ with $x\in T$ and $f_i(x,y)=0$ in $E$.
(Thus $C$ is $\Pi_1^T$, as well as $\Sigma_1$.)  Indeed, for $i\notin C$,
the isomorphic image in $E$ of the elements $(x_i,y_i)$ from $F$ will be such a pair.
For the converse, suppose $i\in C$.  Then $x_i$ and $y_i$ were made algebraic
at some stage in the construction of $F$, and by Theorem \ref{thm:noncovering} with
$k=\Q(y_i:i\in C)\subseteq F$ (which is the subfield containing all elements algebraic
over $\Q$), the function field of the collection $\set{f_j}{j\notin C}$ over $k$
does not contain any nontrivial solution to $f_i=0$.  This function field
is isomorphic to $F$ itself, so we have $C\leq_T T$.
\qed\end{pf}

Theorem \ref{thm:singleton} answers the initial question posed above.
Theorems \ref{thm:upcone} and \ref{thm:edegree} will provide further examples.

\begin{corollary}
\label{cor:singleton}
There exists a computable field $F$ such that no computable field
$E$ isomorphic to $F$ has a computable transcendence basis,
nor even a computable transcendence relation.
\qed\end{corollary}

Our next example shows that the relations $T$ and $I$ can also
\emph{intrinsically compute} a c.e.\ degree $\bfc$, in the sense that
the spectra of $T$ and $I$ can equal the upper cone above $\bfc$
(subject to the restriction that these spectra only contain c.e.\ degrees).
Once again, this parallels a result of Downey and Moses in
\cite{DM91} for linear orders.  In Corollary \ref{cor:Delta2degree} below,
we will generalize this result to all $\Delta^0_2$ degrees $\bfc$,
which is not known (to us) to be possible for adjacency on linear orders.

\begin{theorem}
\label{thm:upcone}
For each computably enumerable Turing degree $\bfc$, there exists
a computable field $F$ for which the spectrum of the transcendence relation
$T$ and of the independence relation $I$ are both
$$ \DS{F}T = \DS{F}I = \set{\bfd\geq\bfc}{\bfd\text{~is a c.e.\ degree}}.$$
\end{theorem}
\begin{pf}
Fix some c.e.\ set $C\in\bfc$, with a computable enumeration $\la C_t\ra_{t\in\omega}$
by finite nested sets $C_t$.  The field $F$ for this degree is the field
$$ F=\Qbar(x_{k})[y_{k}]/(f_k(x_k,y_k)),$$
with $k$ ranging over the set $\Cbar\oplus\omega$.
That is, $F$ has the algebraic closure $\Qbar$ as its ground field,
and contains an algebraically independent set $\set{x_{2i}}{i\in\Cbar}\cup\set{x_{2j+1}}{j\in\omega}$
of elements, along with corresponding elements $y_{2i}$ and $y_{2j+1}$
that ``tag'' the individual $x$-elements by forming solutions to $f_{2i}=0$
or $f_{2j+1}=0$.  (The reason for the odd-indexed elements $x_{2j+1}$
will become clear below:  they will give us the upward closure we desire.)

To see that $F$ has a computable presentation, start building
a computable copy of $\Qbar$, with only finitely many elements added
at each stage.  At stage $s+1$, we add new elements $x_s$ and $y_s$ to the field,
with $f_s(x_s,y_s)=0$, and treat $x_s$ as a transcendental over all previously
existing field elements.  For odd values $s=2j+1$, we simply continue
at each subsequent stage to build the field, with $x_s$ remaining transcendental.
For even $s=2i$, at each subsequent stage $t>s+1$, we check whether
$i\in C_t$.  As long as $i\notin C_t$, we simply add to the field the next element
generated by $x_s$, continuing to treat $x_s$ as transcendental over the
preceding elements.  However, for the first $t$ (if any) with $i\in C_t$,
we switch strategies and make $x_s$ a rational number, finding some way
to do this that is consistent with the finite portion of the atomic diagram
of $F$ that has already been defined.  Of course, this also makes $y_s$
algebraic over $\Q$, though not rational.  This enlarges our presentation
of the ground field $\Qbar$, of course, but since only finitely much of $\Qbar$
had been built so far, it is easy to incorporate $x_s$ and $y_s$ into it
and to continue building $\Qbar$, including them, at each subsequent stage.

Now for any computable field $E\cong F$, with transcendence relation $T$,
we can compute $C$ from $T$.  Indeed, by Theorem \ref{thm:noncovering},
a number $i$ lies in $\Cbar$ if and only if
$E$ contains transcendental elements $x$ and $y$ such that $f_{2i}(x,y)=0$,
so $\Cbar$ is $\Sigma_1^T$, while $C$ is $\Sigma_1$.  Thus
$\DS{F}T$ contains only degrees above $\bfc$, and these must all be
c.e.\ degrees, as the relation $T$ is definable in $F$ by a computable
infinitary $\Pi^0_1$ formula.  The same analysis applies to the independence
relation $I$.

To prove the reverse inclusion, let $\bfd$ be any c.e.\ degree that computes $\bfc$,
and fix some c.e.\ set $D\in\bfd$ with computable enumeration $\la D_s\ra_{s\in\omega}$.
We build a specific computable copy $E$ of $F$
in which $T\equiv_T I\equiv_T D$, by a process quite similar to the above construction
of $F$ itself.  $E$ includes a copy of $\Qbar$, built slowly, with only finitely many
elements added at each stage.  Once again, the even-indexed $x_{2i}$ and $y_{2i}$
are added at stage $2i$ and treated as transcendental until $i$ enters $C$,
at which point $x_{2i}$ becomes rational.  The odd-indexed elements $x_{2j+1}$
and $y_{2j+1}$ are added at stage $2j+1$ and treated as transcendental
until (if ever) we reach a stage $s$ with $j\in D_s$.  If such a stage occurs,
then this $x_{2j+1}$ is made rational at that stage (in the same way as with $x_{2i}$
if $i$ enters $C$), and we adjoin to the field new elements $x'_{2j+1}$ and
$y'_{2j+1}$, again with $x'_{2j+1}$ transcendental over all existing elements
of $F$ and with $f_{2j+1}(x'_{2j+1},y'_{2j+1})=0$.  These new elements
will forever remain transcendental over the ground field $\Qbar$, and the original
$x_{2j+1}$ and $y_{2j+1}$ have now been ``swallowed up'' by $\Qbar$.
Thus the $E$ built here is indeed isomorphic to $F$, and is a computable field.
However, from the transcendence relation $T$ on $F$, we can compute $D$,
since $j\in\Dbar$ if and only if the original $x_{2j+1}$ lies in $T$.  Conversely,
from a $D$-oracle we can decide whether $x_{2j+1}$ will ever be swallowed
up by $\Qbar$ or not, and also (since $C\leq_T D$) whether $x_{2i}$ will
remain transcendental in $E$ or not.  Thus $T\equiv_T D$, and so
$\bfd\in\DS{F}T$.

The same argument also shows that $\bfd\in\DS{F}I$,
since the elements $x_s$ that stay transcendental forever form a
transcendence basis for $E$, from which we can compute
the independence relation.  It should be remarked here, as in
Section \ref{sec:Fermat}, that the first transcendental solution
to $f_k=0$ that one finds in $E$ will only be one of the $6q_i^2$
such solutions, but in enumerating a transcendence basis, it is safe
to choose the first coordinate of the first transcendental solution we find,
and then to ignore all other solutions to the same $f_k$, as their coordinates
are all either in $\Qbar$ or interalgebraic with the coordinate we chose.
In fact, since the automorphism group of $E$ acts transitively on these solutions,
there is nothing to distinguish one such choice from another.
\qed\end{pf}

Our next result suggests that many spectra of transcendence relations can be viewed
as upper cones of enumeration degrees.  To be clear, the spectrum is still a
set of Turing degrees, by definition, but the defining property of the spectrum
may be the ability to enumerate a particular set.  (It remains true
that only c.e.\ degrees may lie in $\DS{F}T$, although other Turing degrees may
enumerate the same set.  So the spectrum will never truly be an upper cone
of $e$-degrees.)

\begin{theorem}
\label{thm:edegree}
Let $S$ be any $\Sigma^0_2$ subset of $\omega$.  Then there
exists a computable field $F$ such that
$$ \DS{F}T=\DS{F}I = \set{\text{c.e.\ degrees~}\bfd}{S\in\Sigma_1^{\bfd}}.$$
That is, $\DS{F}T$ contains exactly those c.e.\ degrees that have
the ability to enumerate $S$.
\end{theorem}
\begin{pf}
Since $S$ is $\Sigma^0_2$, there exists a computable total ``chip function''
$h:\omega\to\omega$ such that $S=\set{n}{h^{-1}(n)\text{~is finite}}$.
The field $F$ we use for this set is the field
$$ F=\Qbar(x_{k})[y_{k}]/(f_i(x_k,y_k)),$$
with $k$ ranging over the set $S\oplus\omega$, much as in Theorem
\ref{thm:upcone} but using the set $S$ itself instead of its complement.

To give a computable presentation of $F$, we start building a copy of
the field $\Qbar(x_k)[y_k]/(f_k(x_k,y_k)$ with $k$ ranging over all of $\omega$,
so that every $x_k$ is initially treated as a transcendental.  For odd $k=2j+1$,
$x_k$ stays transcendental throughout this construction.  For even $k=2i$,
we write $x_{2i,0}=x_{2i}$ for the initial element described above.
At each stage $s+1$, we check whether $h(s)=i$.  If not, then we keep
$x_{2i,s+1}=x_{2i,s}$ and continue
to treat it as a transcendental.  If $h(s)=i$, however, then we suspect
that $i$ might not lie in $S$ (since $h^{-1}(i)$ might turn out to be infinite).
In this case we make the current $x_{2i,s}$ into a rational number, consistently
with the finite portion of the atomic diagram of $F$ built so far, and thus make
$y_{2i,s}$ algebraic.  We then adjoin new elements $x_{2i,s+1}$ and $y_{2i,s+1}$
to $F$, treating $x_{2i,s+1}$ as transcendental and setting $f_{2i}(x_{2i,s+1},y_{2i,s+1})=0$.
We continue building $\Qbar$ as the ground field, now incorporating the old
$x_{2i,s}$ and $y_{2i,s}$ into it, and continuing closing $F$ itself under the field operations,
but always adding only finitely many new elements at each stage.  This completes
the construction, and it is clear that $f_{2i}(X,Y)=0$ will have a solution
by transcendental elements in $F$ just if $h^{-1}(i)$ is finite, which is to say,
just if $i\in S$.  The rest of the construction then makes it clear that the field
we have built is a computable copy of the field $F$ described above.

Given any computable field $E\cong F$, let $T$ be the transcendence relation on $E$.
Then, given a $T$-oracle, we may search in $E$ for a solution to $f_{2i}(X,Y)=0$
using transcendental elements $x$ and $y$.  If we find one, then by the definition
of $F$ we know that $i\in S$.  Conversely, if $i\in S$, then such a solution exists,
and we will eventually find it.  Thus $S$ is c.e.\ relative to the degree $\bfd$ of $T$,
as required.

Conversely, fix any c.e.\ degree $\bfd$ such that $S$ is $\Sigma_1^{\bfd}$, and
fix a c.e. set $D\in\bfd$ and a computable enumeration of it.
Also fix an index $e$ such that $S=W_e^D=\dom{\Phi_e^D}$;
we will use this below to give a computable chip function for $S$,
similar to that used in the original computation of $F$ but specific
to this $D$.
To build a computable copy $E$ of $F$ whose transcendence
relation $T$ satisfies $T\equiv_T D$, we use the strategy from
Theorem \ref{thm:upcone}.  The elements $x_{2i,0}$
and $y_{2i,0}$ are defined and initially treated as transcendentals.
However, at each stage $s+1$, the current $x_{2i,s}$ and $y_{2i,s}$ are made
into algebraic elements and replaced by new elements $x_{2i,s+1}$ and $y_{2i,s+1}$
unless $\Phi_{e,s}^{D_s}(i)\converges$ with some use $u$ such that
$D_{s+1}\res u=D_s\res u$.  This is our new chip function for $S$:
if $i\in S=\dom{\Phi_e^D}$, then there will be some $s_0$ such that
we keep $x_{2i,s_0}$ transcendental at all stages $\geq s_0$;
whereas if $i\notin S$, then for every stage $s+1$ at which
$\Phi_{e,s}^{D_s}(i)\converges$ with a use $u$, there must be some $t>s$
with $D_t\res u\neq D_s\res u$, so that $x_{2i,s}$ will be made
algebraic at stage $t+1$ and replaced by a new $x_{2i,t+1}$.

We also revamp the construction for the odd-indexed elements $x_{2j+1}$ and $y_{2j+1}$,
using exactly the same process as in the proof of Theorem \ref{thm:upcone}.
If we ever reach a stage at which $j$ enters $D$,
then we turn $x_{2j+1}$ into a rational number, consistently with the construction so far,
and adjoin a new transcendental $x_{2j+1}'$ and corresponding $y_{2j+1}'$
with $f_{2j+1}(x'_{2j+1},y'_{2j+1})=0$ in $E$.  This completes the construction of $E$,
which is clearly a computable field and isomorphic to $F$.

Now from an oracle for the
transcendence relation $T$ on $E$, we can determine whether $x_{2j+1}$
is algebraic in $E$ or not, thus deciding whether or not $j\in D$.  Thus $D\leq_T T$.
For the reverse reduction, we claim that with a $D$-oracle we can enumerate
a transcendence basis $B$ for $E$, thus deciding the independence relation $I$ on $E$,
which in turn computes $T$.  This will prove $I\equiv_T T\equiv_T D$ as required.
The $D$-oracle allows us to decide, for each $j$, whether $j\in D$, from which we determine
either that $x_{2j+1}$ lies in $B$ (if $j\notin D$) or that $x_{2j+1}'$ does (if $j\in D$,
in which case we identify $x_{2j+1}'$ by waiting for a stage at which $j$ has entered $D$).
Next, for each $i$, we watch for a stage $s$ at which $\Phi_{e,s}^{D_s}(i)\converges$
with a use $u$ such that $D_s\res u=D\res u$.  The $D$-oracle allows us to check this,
and if we ever find such an $s$, then we enumerate $x_{2i,s+1}$ into our basis, since the
computable enumeration of $D$ will never again change below $u$.  (This is where
our argument would fail if $\bfd$ were a properly $\Delta^0_2$ degree, rather than a c.e.\ degree.
With only a computable approximation to $D$, we could not be sure whether $D_s\res u$
would ever again change, even knowing that $D_s\res u=D\res u$.)  Thus we have
enumerated exactly the set of elements $x_k$ given when we first defined the
isomorphism type of $F$ above, and this set is a transcendence basis for $E$.
\qed\end{pf}

\begin{corollary}
\label{cor:Delta2degree}
Let $\bfc$ be any $\Delta^0_2$ Turing degree.  Then there
exists a computable field $F$ such that
$$ \DS{F}T=\DS{F}I = \set{\text{c.e.\ degrees~}\bfd}{\bfc\leq_T\bfd}.$$
In particular, both $\DS{F}T$ and $\DS{F}I$ can fail to contain a least degree.
\end{corollary}
The condition of not containing a least degree also holds for many
of the spectra given in Theorem \ref{thm:edegree}.  In particular,
if no set $A'$ (with $A$ c.e.) is $1$-equivalent to $S$,
then the spectrum has no least degree, as $\deg{D}\in\DS{F}T$ if and only if
$S\leq_1 D'$.  The proof uses
the Sacks Jump Theorem (see \cite{S63}, or \cite[Thm.\ VII.3.1]{S87}),
to avoid the upper cone above a hypothetical least degree.

\begin{pf}
Fix a set $C\in\bfc$, and apply Theorem \ref{thm:edegree} to the
$\Sigma^0_2$ set $S=C\oplus\Cbar$.  The ability to enumerate $S$
is exactly the ability to compute $C$, so the corollary follows.
(To avoid having a least degree in the spectra, just choose a degree $\bfc\leq\bfz'$
that is not c.e.\ and has no least c.e.\ degree above it.)
\qed\end{pf}

Corollary \ref{cor:Delta2degree} extends Theorem \ref{thm:upcone} to
the $\Delta^0_2$ degrees.  The result can be viewed as an
upper-cone result, but in a somewhat odd way.  If the $\bfc$
in the corollary is c.e., then the corollary merely repeats Theorem
\ref{thm:upcone}.  If $\bfc$ is $\Delta^0_2$ but not c.e., then
the degree spectrum is the restriction of the upper cone above $\bfc$
to the c.e.\ degrees, and therefore does not contain the base degree $\bfc$ itself,
nor any other non-c.e.\ degree $\geq_T\bfc$.

\end{document}